\documentclass[12pt]{article}

\usepackage{amsmath}
\usepackage{amssymb}
\usepackage{amsfonts}
\usepackage{amsthm}
\usepackage{color}
\usepackage{setspace}
\usepackage{graphicx}
\usepackage[english]{babel}\usepackage[latin1]{inputenc}\usepackage{times}\usepackage[T1]{fontenc}
\usepackage{epstopdf}
\usepackage{epsfig}

\newtheorem{thm}{\bf{Theorem}}[section]
\newtheorem{lem}[thm]{\bf{Lemma}}

\newtheorem{df}[thm]{\bf{Definition}}
\newtheorem{cor}[thm]{\bf{Corollary}}
\newtheorem{rem}[thm]{\bf{Remark}}

\newtheorem{prop}[thm]{\bf{Proposition}}
\newtheorem{fact}[thm]{\bf{Fact}}

\newenvironment{rcases}
  {\left.\begin{aligned}}
  {\end{aligned}\right\rbrace}

\numberwithin{equation}{section}

\newcommand{\dom}{\operatorname{dom}}

\newcommand{\gr}{\operatorname{gr}}
\newcommand{\Fix}{\operatorname{Fix}}
\newcommand{\Id}{\operatorname{Id}}

\newcommand{\argmin}{\operatornamewithlimits{argmin}}

\title{Most Convex Functions Have Unique Minimizers\\}
\author{C. Planiden\thanks{Mathematics, University of British Columbia Okanagan, Kelowna, B.C. V1V 1V7, Canada. Research by this author was supported by UBC UGF and by NSERC of Canada. chayne.planiden@lumni.ubc.ca.}\and X. Wang\thanks{Mathematics, University of British Columbia Okanagan, Kelowna, B.C. V1V 1V7, Canada. Research by this author was partially supported by an NSERC Discovery Grant. shawn.wang@ubc.ca.}}
\date{\today}
\begin{document}

\maketitle\author
\setcounter{page}{1}\pagenumbering{arabic}

\begin{abstract}
Finding the minimum and the minimizers of convex functions has been of primary concern in convex analysis since its conception. It is well-known that if a convex function has a minimum, then that minimum is global. The minimizers, however, may not be unique. There are certain subclasses, such as strictly convex functions, that do have unique minimizers when the minimum exists, but other subclasses, such as constant functions, that do not. This paper addresses the question of how many convex functions have unique minimizers. We show, using Baire category theory, that the set of proximal mappings of convex functions that have a unique fixed point is generic.
Consequently, the set of classes of convex functions that have unique minimizers is generic.
\end{abstract}



\section{Introduction}\label{sec:intro}
Convex functions are important in optimization \cite{convanalrock,convmono,rockwets,borweinzhu,vanderwerff}.
This paper builds on the work done in \cite{mostmax}, where it was shown that most monotone operators have unique zeroes. We are concerned with a similar question here: do most convex functions have unique minimizers, or equivalently
do most subdifferential mappings of proper, convex, lsc functions have unique zeros? In fact they do; this is the main result presented in this work. By `most', we are referring to the idea of a generic set in the sense of Baire category. In terms of
 proximal mappings, our result means that
 the set of proximal mappings that do not have unique fixed points is a small, negligible set.
 To prove this, we construct a complete metric space using the proximal mapping as a component of the distance function, and we use an argument based on super-regular mappings  and the density of contraction mappings in the established metric space. Super-regular is a term coined by Reich and Zaslavski \cite[page 2]{genbook}. For a comprehensive study of the
generic properties of nonexpansive mappings, see \cite{genbook}. From the point of view of convex functions, we work with
the set of equivalence classes of functions, as defined in the next section.

\hbox{}
The organization of this paper is as follows. Section \ref{sec:prelim} gives a number of definitions and facts that we will need for the proofs of the main results. Section \ref{sec:prox} contains the definition of the metric space within which we work, and proves that it is complete. Section \ref{sec:main} states and proves the main results of the paper, and Section \ref{sec:conc} makes some concluding remarks and suggests areas of future research.

\section{Preliminaries}\label{sec:prelim}

In this section, we state some definitions and facts that we will use to prove our main results.

\subsection{Notation}
The extended real line $\mathbb{R}\cup\{\infty\}$ is denoted $\overline{\mathbb{R}}.$ The Euclidean space $\mathbb{R}^n$ is
equipped with inner product $\langle x,y\rangle:=\sum\limits_{i=1}^nx_iy_i$ for $x=(x_1,\ldots,x_n),$ $y=(y_1,\ldots,y_n)\in\mathbb{R}^n,$ and induced norm
$\|x\|:=\sqrt{\langle x,x\rangle}.$
We use $\Gamma_0(\mathbb{R}^n)$ to represent the set of proper, convex,
lower semicontinuous (lsc) functions on $\mathbb{R}^n,$ and $B_s(x)$ to represent the closed ball centred at $x$ with radius $s.$ The symbol $\overset{p}\rightarrow$ indicates pointwise convergence. Let $X, Y$ be Hilbert spaces.
We denote the set of continuous functions $f:X\rightarrow Y$ as $C(X,Y).$
The identity function from $X$ to $X$ is denoted $\Id,$ and $\Fix(T)$ is the set of fixed points of the operator $T: X\to X.$ For a sequence $x_{k+1}:=\underset{(k+1\mbox{ times)}}{T\circ T\circ\cdots\circ T}x_0,$ where
 $x_{0}\in X$, we write $x_{k+1}=T^{k+1}(x_0),$ and $g_{k+1}=T^{k+1}.$

\subsection{Definitions}
We start with some concepts about convex functions.
\begin{df}
A function $f\in\Gamma_0(\mathbb{R}^n)$ is \emph{strongly convex} if there exists $\sigma>0$
such that $f-\frac{\sigma}{2}\|\cdot\|^2$ is convex.
\end{df}

\begin{df}\label{superconvex}
A function $f\in\Gamma_0(\mathbb{R}^n)$ is called \emph{super-convex} if it has a unique minimizer.
\end{df}

\begin{df}
A vector $v\in X$ is a \emph{subgradient} of a proper, lsc, convex function $f$ at $\bar{x}$ if
$$f(x)\geq f(\bar{x})+\langle v,x-\bar{x}\rangle\mbox{ for all }x\in\dom f.$$
The \emph{subdifferential $\partial f(\bar{x})$} of a proper, lsc, convex function $f$ at $\bar{x}$ is the set of all subgradients of $f$ at $\bar{x}.$
\end{df}
\begin{df}\label{df:moreau} The \emph{Moreau envelope} of a proper, lsc function $f:X \rightarrow\overline{\mathbb{R}}$ is denoted $e_\lambda f$ and is defined
$$e_\lambda f(x):=\inf\limits_y\left\{f(y)+\frac{1}{2\lambda}\|y-x\|^2\right\}.$$
The parameter $\lambda>0$ is called the \emph{prox-parameter}, and $x$ is called the \emph{prox-center}.
The \emph{proximal mapping} is the (possibly empty) set of points at which this infimum is achieved, and is denoted $P_\lambda f:$
$$P_\lambda f(x):=\argmin\limits_y\left\{f(y)+\frac{1}{2\lambda}\|y-x\|^2\right\}.$$
\end{df}

Next we introduce some terminologies on multi-functions.
\begin{df}
A multifunction $A:X\rightrightarrows X$ is \emph{monotone} if
$$\langle x-y,x^*-y^*\rangle\geq0$$
for all $(x,x^*),(y,y^*)\in\gr A,$ the graph of $A.$ The monotone multifunction $A$ is \emph{maximally monotone} if there does not exist a proper extension of $A$ that is monotone.
\end{df}
\begin{df}
A multifunction $A:X\rightrightarrows X$ is \emph{$n$-cyclically monotone} if
$$\begin{rcases}
(a_i,a_i^*)&\in\gr A,&i=1,2,\ldots,n\\
a_{n+1}&=a_1
\end{rcases}\Rightarrow\sum\limits_{i=1}^n\langle a_{i+1}-a_i,a_i^*\rangle\leq0.$$
When this holds for all $n\in\mathbb{N},$ we say that $A$ is cyclically monotone. We call $A$ maximally cyclically monotone if there does not exist a proper extension of $A$ that is cyclically monotone. It is clear that monotone and 2-cyclically monotone are equivalent.
\end{df}
\begin{df}
The \emph{resolvent} of an operator $A: X\rightrightarrows X$ is
$J_A :=(\Id+A)^{-1}.$
\end{df}

To study uniform convergence on bounded sets of a sequence of mappings we require:
\begin{df}
Let $F\subseteq C(X,Y)$ be nonempty. For a given $x\in X,$ the set $F$ is said to be \emph{equicontinuous at $x$} if to each $\epsilon>0$ there corresponds a neighborhood $U$ of $x$ such that
$$\rho(f(t),f(x))<\epsilon$$
whenever $t\in U$ and $f\in F,$ where $\rho$ is the metric on $Y.$
We say that $F$ is \emph{equicontinuous on $X$} if $F$ is equicontinuous at each $x\in X.$
\end{df}

\begin{df}\label{super}
Let $T:X\rightarrow X$ be a mapping, and $g_n(x)=T^nx.$ The mapping $T$ is called \emph{super-regular} if there exists a unique $x_T\in X$ such that for each $s>0,$ when $n\rightarrow\infty,$ the sequence $\{g_n\}_{n=1}^\infty$ converges to the constant function $x_T$ uniformly on $B_s(0).$
\end{df}

The key concept we need is the Baire category which comes as follows.
\begin{df}\label{nowheredense}
A set $S\subseteq X$ is \emph{nowhere dense} if the interior of its closure is empty.
A set $S\subseteq X$  is \emph{of first category (meagre)} if $S$ is a union of countably many nowhere dense sets.
A set $S\subseteq X$ is \emph{generic or residual} if $X\setminus S$ is of first category.
\end{df}

\begin{df}\label{bairespace}
A topological space $S\subseteq X$ is called a \emph{Baire space} if every intersection of countably many dense, open sets in $S$
is dense in $S$.
\end{df}

\subsection{Facts}
Let us collect some well-known facts used in later proofs.
Properties on convex functions come first.
\begin{fact}\cite[Corollary 3.37]{rockwets}\label{factprox}
If $f_1,$ $f_2\in\Gamma_0(\mathbb{R}^n)$ with $P_\lambda f_1=P_\lambda f_2$ for some $\lambda>0,$ then $f_1=f_2+c,$ where $c
\in\mathbb{R}$ is a constant.
\end{fact}

\begin{fact}\cite[Proposition 12.19]{rockwets}\label{factpmaxmono}
If a function $f:\mathbb{R}^n\rightarrow\overline{\mathbb{R}}$ is proper, lsc, and convex, then $P_\lambda f$ is maximally monotone and non-expansive for all $\lambda>0.$ Hence, $P_\lambda f$ is single-valued.
\end{fact}

\begin{fact}\cite[Theorem 12.17]{rockwets}\label{factsubdiff}
A function $f:\mathbb{R}^n\rightarrow\overline{\mathbb{R}}$ is proper, lsc and convex if and only if $\partial f:\mathbb{R}^n\rightarrow\mathbb{R}^n$ is monotone, in which case $\partial f$ is maximally monotone.
\end{fact}

\begin{fact}\cite[Theorem 12.25]{rockwets}\label{facttcyc}
An operator $T$ is the subdifferential of some $f\in\Gamma_0(\mathbb{R}^n),$ if and only if $T$ is maximally cyclically monotone. Then $f$ is uniquely determined by $T,$ up to a constant.
\end{fact}

\begin{fact}\cite[Example 23.3]{convmono}\label{factrespartial}
If $f:X\rightarrow\overline{\mathbb{R}}$ is proper, lsc, and convex, and $\lambda>0,$ then
$$J_{\lambda\partial f}=P_\lambda f.$$
\end{fact}

Properties on resolvents of monotone operators are:
\begin{fact}\cite[Theorem 6.6]{fitz}\label{factresolv}
Let $X$ be a real Hilbert space, $T:X\rightarrow X.$ Then $T$ is the resolvent of the maximally cyclically monotone operator $A:X\rightrightarrows X$ if and only if $T$ has full domain, $T$ is firmly nonexpansive, and for every set of points $\{x_1,x_2,\ldots,x_n\}$ where $n>1$ and $x_{n+1}=x_1,$ one has
$$\sum\limits_{i=1}^n\langle x_i-Tx_i,Tx_i-Tx_{i+1}\rangle\geq0.$$
\end{fact}

\begin{fact}\cite[Fact 6.2]{mostmax}\label{factnonexp}
Suppose that $X$ is a real Hilbert space, $C\subset X,$ and $T:C\rightarrow X.$ Then the following are equivalent.
\begin{itemize}
\item[i)] $T$ is firmly nonexpansive.
\item[ii)] $\|Tx-Ty\|^2\leq\langle x-y,Tx-Ty\rangle$ for all $x,y\in C.$
\item[iii)] $T=\frac{1}{2}\Id+\frac{1}{2}N,$ where $N$ is nonexpansive (i.e. 1-Lipschitz continuous).
\item[iv)] $T=(\Id+A)^{-1}$ is the resolvent of some monotone multifunction $A:X\rightrightarrows X.$
\end{itemize}
\end{fact}

\begin{fact}\cite[Proposition 1.5]{mostmax}\label{factmost}
Let $X$ be a real Hilbert space, and $M(X)$ be the set of maximally monotone operators on $X.$ The following are equivalent:
\begin{itemize}
\item[a)] a sequence of maximally monotone operators $(A_n)_{n=1}^\infty\in M(X)$ converges graphically to $A,$
\item[b)] $(J_{A_n})_{n=1}^\infty$ converges pointwise to $J_A$ on $X.$
\end{itemize}
\end{fact}

Uniform convergence of a sequence of mappings or functions are given by
\begin{fact}\cite[Theorem 3.1]{genbook}\label{fact:ss}
Let $K$ be a bounded, closed, convex subset of $X,$ and $\mathcal{A}$ be the set of all operators $A:K\rightarrow K$ such that
$$\|Ax-Ay\|\leq\|x-y\|\mbox{ for all }x,y\in K.$$
Assume that $B\in\mathcal{A}$ is a contractive mapping. Then there exists $x_B\in K$ such that $B^nx\rightarrow x_B$ as $n\rightarrow\infty,$ uniformly on $K.$
\end{fact}

\begin{fact}\cite[Theorem 3.143]{introclass}\label{factequi}
Suppose that the metric space $Y$ is complete and that $\{f_n\}_{n=1}^\infty$ is an equicontinuous sequence in $C(X,Y)$ that converges at each point of a dense subset $D$ of the topological space $X.$ Then there is a function $f\in C(X,Y)$ such that $f_n\rightarrow f$ uniformly on each compact subset of $X.$
\end{fact}

The following is crucial for our Baire category techniques in this paper.

\begin{fact}\cite[Theorem 10.11.4]{metricspaces}
Every complete metric space is a Baire space.
\end{fact}

\section{The Complete Metric Space of Subdifferentials}\label{sec:prox}

In this section we establish a complete metric space whose distance function makes use of the proximal mapping.
Many nice properties about proximal mappings can be found in Moreau's seminarial paper \cite{moreau},
as well as that of Combettes and Wajs \cite{combettes}.
In order to prove completeness, we first state the following lemma.
\begin{lem}\label{lem:m4}
Define $a:[0,\infty)\rightarrow\mathbb{R},~a(t):=\frac{t}{1+t}.$ Then
\begin{itemize}
\item[a)] $a$ is an increasing function, and
\item[b)] $t_1,t_2\geq0\Rightarrow a(t_1+t_2)\leq a(t_1)+a(t_2).$
\end{itemize}
\end{lem}
\textbf{Proof.}
\begin{itemize}
\item[a)] $a'(t)=(1+t)^{-2}>0~\forall t\in[0,\infty),$ hence, $a$ is increasing everywhere.
\item[b)] $a(t_1+t_2)=\frac{t_1+t_2}{1+(t_1+t_2)}=\frac{t_1}{1+(t_1+t_2)}+\frac{t_2}{1+(t_1+t_2)}\leq\frac{t_1}{1+t_1}+\frac{t_2}{1+t_2}=a(t_1)+a(t_2).$\qed
\end{itemize}

Let $(J,d)$ be a space defined by
\begin{align*}
J&:=\{\partial f:f\in\Gamma_0(\mathbb{R}^n)\},\mbox{ and}\\
d(\partial f,\partial g)&:=\sum\limits_{i=1}^\infty\frac{1}{2^i}\frac{\sup\limits_{\|x\|\leq i}|P_1f(x)-P_1g(x)|}{1+\sup\limits_{\|x\|\leq i}|P_1f(x)-P_1g(x)|}.
\end{align*}
\begin{prop}
$(J,d)$ is a complete metric space.
\end{prop}
\textbf{Proof.} Items M1-M4 show that $(J,d)$ is a metric space, and item C shows that it is complete.
\begin{itemize}
\item[M1:] \begin{align*}
\sum\limits_{i=1}^\infty\frac{1}{2^i}&=1,\mbox{ and }0\leq\frac{\sup\limits_{\|x\|\leq i}|P_1f(x)-P_1g(x)|}{1+\sup\limits_{\|x\|\leq i}|P_1f(x)-P_1g(x)|}<1,\\
\Rightarrow\frac{1}{2^i}&\geq\frac{1}{2^i}\frac{\sup\limits_{\|x\|\leq i}|P_1f(x)-P_1g(x)|}{1+\sup\limits_{\|x\|\leq i}|P_1f(x)-P_1g(x)|}~\forall i\\
\Rightarrow0&\leq d(\partial f,\partial g)\leq1~\forall\partial f,\partial g\in J.
\end{align*}
Hence, $d$ is real-valued, finite, and non-negative.
\item[M2:] \begin{align*}
d(\partial f,\partial g)=0&\Leftrightarrow\sum\limits_{i=1}^\infty\frac{1}{2^i}\frac{\sup\limits_{\|x\|\leq i}|P_1f(x)-P_1g(x)|}{1+\sup\limits_{\|x\|\leq i}|P_1f(x)-P_1g(x)|}=0\\
&\Leftrightarrow\sup\limits_{\|x\|\leq i}|P_1f(x)-P_1g(x)|=0~\forall i\\
&\Leftrightarrow P_1f(x)-P_1g(x)=0~\forall x\\
&\Leftrightarrow f=g+c\mbox{ (Fact \ref{factprox}).}\\
&\Leftrightarrow\partial f=\partial g.
\end{align*}
Hence $ d(f,g)=0$ if and only if $\partial f=\partial g.$
\item[M3:] $d(\partial f,\partial g)=d(\partial g,\partial f)$ is trivial.
\item[M4:] Let $\partial f,\partial g,\partial h\in J.$ By the triangle inequality for real numbers, we know that
$$|P_1f(x)-P_1g(x)|\leq|P_1f(x)-P_1h(x)|+|P_1h(x)-P_1g(x)|~\forall\partial f,\partial g,\partial h\in J.$$
This gives
\begin{align*}
\sup\limits_{\|x\|\leq i}|P_1f(x)-P_1g(x)|&\leq\sup\limits_{\|x\|\leq i}(|P_1f(x)-P_1h(x)|+|P_1h(x)-P_1g(x)|)\\
&\leq\sup\limits_{\|x\|\leq i}|P_1f(x)-P_1h(x)|+\sup\limits_{\|x\|\leq i}|P_1h(x)-P_1g(x)|.
\end{align*}
By applying Lemma \ref{lem:m4}(a), we see that
$$\frac{\sup\limits_{\|x\|\leq i}|P_1f(x)-P_1g(x)|}{1+\sup\limits_{\|x\|\leq i}|P_1f(x)-P_1g(x)|}\leq\frac{\sup\limits_{\|x\|\leq i}|P_1f(x)-P_1h(x)|+\sup\limits_{\|x\|\leq i}|P_1h(x)-P_1g(x)|}{1+\sup\limits_{\|x\|\leq i}|P_1f(x)-P_1h(x)|+\sup\limits_{\|x\|\leq i}|P_1h(x)-P_1g(x)|}.$$
Then we use Lemma \ref{lem:m4}(b), with $t_1=\sup\limits_{\|x\|\leq i}|P_1f(x)-P_1h(x)|$ and $t_2=\sup\limits_{\|x\|\leq i}|P_1h(x)-P_1g(x)|,$ and we have
$$\frac{\sup\limits_{\|x\|\leq i}|P_1f(x)-P_1g(x)|}{1+\sup\limits_{\|x\|\leq i}|P_1f(x)-P_1g(x)|}\leq\frac{\sup\limits_{\|x\|\leq i}|P_1f(x)-P_1h(x)|}{1+\sup\limits_{\|x\|\leq i}|P_1f(x)-P_1h(x)|}+\frac{\sup\limits_{\|x\|\leq i}|P_1h(x)-P_1g(x)|}{1+\sup\limits_{\|x\|\leq i}|P_1h(x)-P_1g(x)|}.$$
Multiplying both sides by $\frac{1}{2^i}$ and taking the infinite summation over $i$ of both sides, we obtain the distance functions, which yields  $d(\partial f,\partial g)\leq d(\partial f,\partial h)+d(\partial h,\partial g)~\forall \partial f,\partial g,\partial h\in J.$

Combining M1-M4, we see that $(J,d)$ is a metric space.
\item[C:] Let $(\partial f_k)$ be a Cauchy sequence in $(J,d).$ Then $\forall\epsilon>0~\exists N_\epsilon\in\mathbb{N}$ such that $d(\partial f_j,\partial f_k)<\epsilon~\forall j,k\geq N_\epsilon.$ Fix $\epsilon>0.$ Then $\exists N\in\mathbb{N}$ such that
$$\sum\limits_{i=1}^\infty\frac{1}{2^i}\frac{\sup\limits_{\|x\|\leq i}|P_1f_j(x)-P_1f_k(x)|}{1+\sup\limits_{\|x\|\leq i}|P_1f_j(x)-P_1f_k(x)|}<\epsilon~\forall j,k\geq N.$$
Then for $i$ fixed, we have $\frac{1}{2^i}\frac{\sup\limits_{\|x\|\leq i}|P_1f_j(x)-P_1f_k(x)|}{1+\sup\limits_{\|x\|\leq i}|P_1f_j(x)-P_1f_k(x)|}<\epsilon,$ so that $\sup\limits_{\|x\|\leq i}|P_1f_j(x)-P_1f_k(x)|<\frac{2^i\epsilon}{1-2^i\epsilon}.$ Since $\epsilon>0$ is arbitrary, we have that $(P_1f_k(x))_{k=1}^\infty$ is a Cauchy sequence on $\|x\|\leq i,$ so that $P_1f_k(x)\overset{p}\rightarrow h$ for some $h.$ By Fact \ref{factsubdiff}, $\partial f_k$ is maximally monotone for all $k.$ Since $P_1f_k$ is the resolvent of the maximally cyclically monotone operator $\partial f_k,$ the domain of $P_1f_k$ is $\mathbb{R}^n.$ By Fact \ref{factnonexp} we have
\begin{equation}\label{eq:nonexp}
\|P_1f_k(x)-P_1f_k(y)\|^2\leq\langle x-y,P_1f_k(x)-P_1f_k(y)\rangle\mbox{ for all }x,y\in\mathbb{R}^n.
\end{equation}
Then letting $k\rightarrow\infty,$ we get
$$\|h(x)-h(y)\|^2\leq\langle x-y,h(x)-h(y)\rangle\mbox{ for all }x,y\in\mathbb{R}^n,$$
so $h=J_A$ for an maximally monotone operator $A:\mathbb{R}^{n}\rightrightarrows\mathbb{R}^n$.
It remains to be shown that $h=J_{\partial f}$ for some convex function $f\in\Gamma_{0}(\mathbb{R}^n).$
By Fact \ref{factresolv} we have
\begin{equation}\label{eq:nonexp2}
\sum\limits_{i=1}^n\langle x_i-P_1f_k(x_i),P_1f_k(x_i)-P_1f_k(x_{i-1})\rangle\geq0
\end{equation}
for all $\{x_1,x_2,\ldots,x_n\}$ with $x_{n+1}=x_1,$ for all $n>1.$ Then, letting $k\rightarrow\infty$ in equations (\ref{eq:nonexp}) and (\ref{eq:nonexp2}), we get
\begin{equation}\label{eq:nonexp3}
\|h(x)-h(y)\|^2\leq\langle x-y,h(x)-h(y)\rangle
\end{equation}
for all $x,y\in\mathbb{R}^n,$ and
\begin{equation}\label{eq:nonexp4}
\sum\limits_{i=1}^n\langle x_i-h(x_i),h(x_i)-h(x_{i-1})\rangle\geq0
\end{equation}
for all $\{x_1,x_2,\ldots,x_n\}$ with $x_{n+1}=x_1,$ for all $n>1.$ Hence, $h$ is a cyclical resolvent.
Equation (\ref{eq:nonexp3}) says that $P_1f$ is full-domain and firmly nonexpansive, which together with equation (\ref{eq:nonexp4}) gives us via Fact \ref{factresolv} that $h=J_A$ is the resolvent of a maximally cyclically monotone operator. Hence, $A$ is maximally cyclically monotone, which means that $A$ is the subdifferential of a proper, lsc, convex function by Fact \ref{facttcyc}. Therefore, $(J,d)$ is closed, and is a complete metric space.\qed
\end{itemize}

\begin{rem}
By Fact \ref{factmost} in $(J,d),$ the convergence is graphical convergence. Thus,
the topological space $(J,\Gamma)$ where $\Gamma$ denotes graphical convergence, is metrizable.
\end{rem}
We proceed to introduce two closely related metric spaces.

For all $f\in\Gamma_0(\mathbb{R}^n),$ define the equivalence classes $\mathbb{F}_f:$
$$\mathbb{F}_f:=\{g\in\Gamma_0(\mathbb{R}^n):f-g=c, \text{where $c\in\mathbb{R}$ is a constant}\}.$$
 We denote by $\mathcal{F}$ the set of all such equivalence classes:
$$\mathcal{F}:=\{\mathbb{F}_f:f\in\Gamma_0(\mathbb{R}^n)\}.$$
This forms a partition of $\Gamma_0(\mathbb{R}^n).$ That is,
the intersection of any two distinct elements of $\mathcal{F}$ is empty, and $\Gamma_0(\mathbb{R}^n)=\bigcup\limits_{\mathbb{F}_f\in\mathcal{F}}\mathbb{F}_f.$
Now considering the equivalence classes $\mathbb{F}_f$ and $\mathbb{F}_g,$ define the metric $\tilde{d}:$
$$\tilde{d}(\mathbb{F}_f,\mathbb{F}_g):=\sum\limits_{i=1}^\infty\frac{1}{2^i}\frac{\sup\limits_{\|x\|\leq i}\|P_1f(x)-P_1g(x)\|}{1+\sup\limits_{\|x\|\leq i}\|P_1f(x)-P_1g(x)\|},$$
where $f$ and $g$ are arbitrary elements of $\mathbb{F}_f, \mathbb{F}_g \in\mathcal{F}$ respectively.
Then by Fact~\ref{factprox}, we have that $\tilde{d}$ is a metric on $\mathcal{F},$ and
$$\tilde{d}(\mathbb{F}_f,\mathbb{F}_g)=d(\partial f,\partial g).$$
Thus the following corollary holds.

\begin{cor}
The space $(\mathcal{F},\tilde{d})$ is a complete metric space.
\end{cor}

Define
$\mathcal{P}:=\{P_1f:f\in\Gamma_0(\mathbb{R}^n)\},$ and
$$\rho(T_{1},T_{2}):=
\sum\limits_{i=1}^\infty\frac{1}{2^i}\frac{\sup\limits_{\|x\|\leq i}\|T_{1}x-T_{2}x\|}{1+
\sup\limits_{\|x\|\leq i}\|T_{1}x-T_{2}x\|}$$
where
$T_{1},T_{2}\in\mathcal{P}$.
\begin{cor}
The space $(\mathcal{P},\rho)$ is a complete metric space.
\end{cor}

Although $(J,d), (\mathcal{F},\tilde{d}), (\mathcal{P},\rho)$ look different, they are in fact isometric.
\begin{prop}
The three complete metric spaces $(J, d)$, $(\mathcal{F}, \tilde{d})$ and $(\mathcal{P}, \rho)$
are isometric.
\end{prop}

\textbf{Proof.} Define $\phi:(J,d)\rightarrow (\mathcal{P},\rho)$ by
$$\phi(\partial f)=P_1f=(\partial f+\Id)^{-1}.$$
Then $\rho(\phi(\partial f),\phi(\partial g))=d(\partial f,\partial g)$ for all $\partial f, \partial g\in J$,
 and $\phi$ is bijective. Therefore, $(J,d)$ and $(\mathcal{P},\rho)$ are isometric.

Define $\psi:(\mathcal{F},\tilde{d})\rightarrow (J, d)$ by
$$\psi(\mathbb{F}_f)=\partial f.$$
Then $d(\psi(\mathbb{F}_f),\psi(\mathbb{F}_g))=\tilde{d}(\mathbb{F}_f,\mathbb{F}_g)$ for all
$\mathbb{F}_{f},\mathbb{F}_{g}\in\mathcal{F}$, and $\psi$ is bijective
by Fact~\ref{facttcyc}. Therefore, $(\mathcal{F},\tilde{d})$ and $ (J, d)$ are isometric.
\qed

\section{Main Results}\label{sec:main}

In this section we establish the main result: the set of proximal mappings that have a unique fixed point is a generic set. Equivalently, the set of equivalence classes of convex functions that have unique minimizers is a generic set. For all that follows, we use $f$ to represent any function in $\mathbb{F}_f,$ as the results are the same for any function in the equivalence class of $f.$ To start, we need the following results, which give conditions for super-regularity of the proximal mapping and show that the set of contractive proximal mappings is dense in $\mathcal{P}=\{P_1f:f\in\Gamma_0(\mathbb{R}^n)\}.$ We will use these in the proof of the main result.
\begin{prop}\label{prop:strictsuper}
If $f\in \Gamma_{0}(\mathbb{R}^n)$ is strictly convex and $\argmin f\neq\emptyset,$ then $P_1f$ is super-regular.
\end{prop}
\textbf{Proof.} Since $f$ is strictly convex and $\argmin f\neq\emptyset,$ then the minimizer is a singleton, $\argmin f=\{x_T\}.$ By the \emph{proximal point algorithm} \cite[Theorem 4]{monops}, we know that by starting at any arbitrary point $x$ and iteratively calculating $x_{k+1}=(P_1f)^k(x)=g_k(x),$ this generates a sequence of functions $\{g_k(x)\}$ that converges to $x_T,$ where $g_k(x)=P_1f(x_k).$ Since $\{g_k(x)\}$ is a collection of non-expansive mappings, they are Lipschitz-1 and hence equicontinuous on $\mathbb{R}^n$ \cite[Theorem 12.32]{rockwets}. So by Fact \ref{factequi}, we have that $(P_1f)^k(x)\rightarrow x_T$ uniformly on each compact subset of $\mathbb{R}^n.$ In particular, $(P_1f)^k\rightarrow x_T$ uniformly on $B_s(0)$ for all $s>0.$ Hence, $P_1f$ is super-regular.\qed\\

Proposition \ref{prop:strictsuper} is a particular case (strictly convex functions), but it offers us a hint at the direction in which to go to obtain a necessary and sufficient condition for super-regular proximal mappings. That condition is found in Theorem \ref{thm:super} below.
\begin{thm}\label{thm:super}
Let $f\in\Gamma_{0}(\mathbb{R}^n).$ Then $P_1f$ is super-regular if and only if $\argmin f$ is a singleton.
\end{thm}
\textbf{Proof.}
\begin{itemize}
\item[$(\Rightarrow)$] Suppose $P_1f$ is super-regular. Define $\{g_k(x)\}$ as the iterative sequence $g_k(x)=(P_1f)^k(x)$ as in Proposition \ref{prop:strictsuper}. Then there exists a unique $x_T\in\mathbb{R}^n$ such that $g_k(x)\rightarrow x_T$ uniformly on $B_s(0)$ for any $s>0.$ By the proximal point algorithm, we know that $\{g_k(x)\}$ converges to the minimizer
    of $f.$ Therefore, $\argmin f=\{x_T\},$ a singleton.
\item[$(\Leftarrow)$] Suppose $\argmin f=\{x_T\}$ is a singleton. By the proximal point algorithm, $g_k(x)\rightarrow x_T.$ In the proof of Proposition \ref{prop:strictsuper}, we saw via equicontinuity that this convergence is uniform on each compact subset of $\mathbb{R}^n,$ and in particular on $B_s(0)$ for all $s>0.$ Therefore, $P_1f$ is super-regular.\qed
\end{itemize}
\textbf{Remark:} The set of super-regular proximal mappings is strictly larger than the set of proximal mappings of strictly convex functions with unique minimizers. For example, the function $f(x)=|x|$ has a super-regular proximal mapping, yet $f(x)$ is not strictly convex.
\begin{lem}\label{lem:prox1}
In $(J,d),$ the set of strongly monotone mappings is dense. Equivalently, in $(\mathcal{F},\tilde{d})$, the
set of strongly convex function classes is dense; in $(\mathcal{P},\rho)$, the set of contraction mappings
is dense.
\end{lem}
\textbf{Proof.} We need to show that for every $\epsilon>0$ and $T\in\mathcal{P},$ there exists a contraction $\overline{T}$ such that $d(\overline{T},T)<\epsilon.$ As $T\in\mathcal{P},$ $T=P_1f$ for some $f\in\Gamma_0(X).$ Define $\overline{T}:=(1-\sigma)P_1f$ for some $f\in\Gamma_0(\mathbb{R}^n),$ $\sigma\in(0,1).$ Then $\overline{T}$ is a contraction. Our first goal is to find a function $g\in\Gamma_0(\mathbb{R}^n)$ such that $\overline{T}=P_1g.$ We do this by equating $\overline{T}$ to the resolvent of $g,$ and solving for $g.$ This follows from:
\begin{align*}
P_1g&=(1-\sigma)P_1f\\
(\Id+\partial g)^{-1}&=(1-\sigma)(\Id+\partial f)^{-1}\\
(\Id+\partial g)^{-1}&=\left[(\Id+\partial f)\circ\left(\frac{1}{1-\sigma}\Id\right)\right]^{-1}\\
(\Id+\partial g)(x)&=(\Id+\partial f)\left(\frac{1}{1-\sigma}\Id\right)(x)\\
x+\partial g(x)&=\frac{x}{1-\sigma}+\partial f\left(\frac{x}{1-\sigma}\right)\\
\partial g(x)&=\frac{\sigma}{1-\sigma}x+\partial f\left(\frac{x}{1-\sigma}\right)\\
\partial g&=\frac{\sigma}{1-\sigma}\Id+\partial f\circ\left(\frac{1}{1-\sigma}\Id\right).
\end{align*}
From here we see that $\partial g$ is strongly monotone, so that $g\in\Gamma_0(\mathbb{R}^n)$
 is strongly convex. Thus, $\overline{T}$ is the proximal mapping of the proper, lsc, strongly convex function $g,$ where
$$g(x)=\frac{\sigma}{1-\sigma}\frac{\|x\|^2}{2}+(1-\sigma)f\left(\frac{x}{1-\sigma}\right).$$
For $\epsilon>0,$ choose $N$ such that $\sum\limits_{i\geq N}^\infty\frac{1}{2^i}<\frac{\epsilon}{2}.$ Consider
\begin{align*}
d(\partial f, \partial g)= \rho(P_1f,P_1g)&\leq\sum\limits_{i=1}^N\frac{1}{2^i}\frac{\sup\limits_{\|x\|\leq i}\|P_1f(x)-P_1g(x)\|}{1+\sup\limits_{\|x\|\leq i}\|P_1f(x)-P_1g(x)\|}+\frac{\epsilon}{2}.
\end{align*}
Notice that
\begin{align*}
\sup\limits_{\|x\|\leq i}\|P_1f(x)-P_1g(x)\|&=\sup\limits_{\|x\|\leq i}\|(1-\sigma)P_1f(x)-P_1f(x)\|\\
&=\sigma\sup\limits_{\|x\|\leq i}\|P_1f(x)\|.
\end{align*}
As $P_1f:\mathbb{R}^n\rightarrow\mathbb{R}^n$ is nonexpansive, there exists $M>0$ such that
$$\sup\limits_{\|x\|\leq i}\|P_1f(x)\|<M.$$
Hence,
$$\sup\limits_{\|x\|\leq i}\|P_1f(x)-P_1g(x)\|<\sigma M.$$
Then
\begin{align*}
\rho(P_1f,P_1g)&\leq\sum\limits_{i=1}^N\frac{1}{2^i}\frac{\sigma M}{1+\sigma M}+\frac{\epsilon}{2}\\
&=\frac{\sigma M}{1+\sigma M}\sum\limits_{i=1}^N\frac{1}{2^i}+\frac{\epsilon}{2}\\
&\leq\frac{\sigma M}{1+\sigma M}+\frac{\epsilon}{2}\\
&\leq\sigma M+\frac{\epsilon}{2}\\
&<\epsilon\mbox{ \qquad\qquad\qquad when }\sigma<\frac{\epsilon}{2M}.
\end{align*}
Thus, for any $\epsilon>0$ one can always choose $\sigma$ small enough so that $d(\partial f,\partial g)<\epsilon$. Therefore,
the set of strongly monotone subdifferential mappings is dense in $(J,d).$ Since $(J,d),$ $(\mathcal{F},\tilde{d}),$ and $(\mathcal{P},\rho)$ are isometric, we have the equivalent conclusions:
\begin{itemize}
\item[i)] in $(\mathcal{F},\tilde{d})$ the set of strongly convex function classes is dense, and
\item[ii)] in $(\mathcal{P},\rho)$ the set of contraction mappings is dense.
\end{itemize}\qed\\
\begin{lem}\label{factsuperregular}
In $(\mathcal{P},\rho)$,
let $T\in\mathcal{P}$ be super-regular,
$\Fix(T)=\{x_T\}.$ Then for every $\epsilon>0$ and $s>0,$ there exist $\delta>0$ and $n_0>1$ such
that when $\rho(\overline{T},T)<\delta$ and $n\geq n_0,$ we have
$$\|\overline{T}^nx-x_T\|<\epsilon\mbox{ for every }x\in B_s(0).$$
\end{lem}
\textbf{Proof.}
Apply \cite[Proposition 2.12]{mostmax} with $\mathcal{F}=\mathcal{P}$ and $d=\rho$.
\qed

The next lemma we need in order to state the main result proves that the set of super-regular proximity operators is a generic set in $\mathcal{P}.$
\begin{lem}\label{lem:prox2}
In $(\mathcal{P},\rho),$ there exists a set $G\subset\mathcal{P}$ that is a countable intersection of open, everywhere dense sets in $\mathcal{P}$ such that
each $T\in G$ is super-regular. Hence, $G$ is a generic subset of $\mathcal{P}.$
\end{lem}
\textbf{Proof.} The proof is similar to that of \cite[Proposition 2.14]{mostmax}. Let $\mathcal{C}$ be the set of contractive
 proximal mappings; recall that $\mathcal{C}$ is dense in $\mathcal{P}$ by Lemma \ref{lem:prox1}. By Fact \ref{fact:ss} or \cite[Proposition 2.11(i)]{mostmax}, each $T\in
\mathcal{C}$ is super-regular. By Lemma~\ref{factsuperregular}, for each $T\in\mathcal{C}$ there exists an open neighborhood $U(T,i)$ of $T$ in $(\mathcal{P},d),$ and a natural number $n(T,i)$ such that
\begin{equation}\label{eq:sup1}
\|\overline{T}^nx-x_T\|<\frac{1}{i}
\end{equation}
whenever $\overline{T}\in U(T,i),$ $n\geq n(T,i),$ and $x\in B_i(0).$ Define
$$O_q:=\bigcup\{U(T,i):T\in\mathcal{C},i\geq q\}.$$
Then $\mathcal{C}\subset O_q.$ Now we define $G:=\bigcap\limits_{q=1}^\infty O_q.$ Since each $U(T,i)$ is open and $\mathcal{C}\subset O_q,$ we have that $G$ is dense in $\mathcal{P}.$\\
It remains to show that every $T\in G$ is super-regular. Let $T\in G$ be arbitrary. Then there exist sequences $\{T_q\}_{q=1}^\infty$ and $\{i_q\}_{q=1}^\infty$ such that  $T\in U(T_q,i_q)$ for each $q\in\mathbb{N}.$ By using inequality (\ref{eq:sup1}), we get
\begin{equation}\label{eq:sup2}
\|T^nx-x_{T_q}\|<\frac{1}{i_q}
\end{equation}
whenever $n\geq n(T_q,i_q)$ and $x\in B_{i_q}(0).$ Hence, letting $N=\max\{n(T_q,i_q),n(T_p,i_p)\}$ and $M=\min\{i_q,i_p\},$ we know that
$$\|x_{T_q}-x_{T_p}\|\leq\|x_{T_q}-T^nx\|+\|T^nx-x_{T_p}\|<\frac{1}{i_q}+\frac{1}{i_p}$$
whenever $n\geq N$ and $\|x\|\leq M.$ In other words, we have a Cauchy sequence $\{x_{T_q}\}_{q=1}^\infty$ such that $x_{T_q}\rightarrow x_T.$ With the correct choice of $q$ and $i_q,$ we are sure that $B_s(0)\subset B_{i_q}(0),$ and that
$$\frac{1}{i_q}+\|x_{T_q}-x_T\|<\epsilon.$$
Now using that fact together with inequality (\ref{eq:sup2}), we have
$$\|T^nx-x_T\|\leq\|T^nx-x_{T_q}\|+\|x_{T_q}-x_T\|<\frac{1}{i_q}+\|x_{T_q}-x_T\|<\epsilon$$
when $n\geq N$ and $\|x\|\leq s.$
Hence, $T$ is super-regular, and since $T$ is an arbitrary element of $G,$ every element of $G$ is super-regular. Therefore, the set of super-regular proximity mappings in $\mathcal{P}$ is a generic set.\qed

We are ready to present the main results.
\begin{thm}\label{main0}
 In $(\mathcal{P},\rho)$, define the set of proximal mappings
$$S:=\{T\in\mathcal{P}: \Fix T \text{ is a singleton}\}.$$ Then
$S$ is generic in $\mathcal{P}$.
\end{thm}
\textbf{Proof.}  Every super-regular mapping $T$ has $\Fix(T)$ a singleton. Since the set
$G$ given in Lemma~\ref{lem:prox2}
satisfies $G\subseteq S,$ we have that $S$ is generic in $\mathcal{P}.$\qed

\begin{thm}\label{main1}
On $(J,d),$ define the set of subdifferentials
$$S:=\{\partial f\in J: f\in \Gamma_{0}(\mathbb{R}^n), \partial f\mbox{ has a unique zero}\}.$$ Then $S$ is generic in $J.$
\end{thm}
\textbf{Proof.} Since every element of $S$ has a unique zero, by Theorem \ref{thm:super} we have that $P_1f$ is super-regular for every corresponding $\partial f\in S.$ Then by Theorem~\ref{main0}, the set $\{P_1f:\partial f\in S\}$ is generic in $\mathcal{P}.$ Since $(J,d)$ and $(\mathcal{P},\rho)$ are isometric, $S$ is generic in $J.$\qed

Also because $(J,d)$ and $(\mathcal{F},\tilde{d})$ are isometric, we obtain:
\begin{thm}
In $(\mathcal{F},\tilde{d}),$ define the set of classes of convex functions
$$S:=\{\mathcal{F}_f: f\in \Gamma_{0}(\mathbb{R}^n), f\mbox{ has a unique minimizer}\}.$$ Then $S$ is generic in $\mathcal{F}.$
\end{thm}

\section{Conclusion}\label{sec:conc}

We have established a complete metric space using proximal mappings, and used it to show that for proper, convex, lsc functions, the set of proximal mappings that do not have a unique fixed point is a Baire category one set.
There are at least two areas of further research to be done in this vein.
\begin{itemize}
\item[1)] All results presented here are for functions on $\mathbb{R}^n.$ With some additional arguments, extension to infinite-dimensional space should be attainable.
\item[2)] This paper was presented from the proximal mapping point of view, but one should be able to reach a similar conclusion by working from the function point of view as well, that is, proving that the set of proper, convex, lsc functions that have a unique minimizer is generic.
\end{itemize}

\bibliographystyle{plain}
\bibliography{Convex}{}
\end{document}